\documentclass[twocolumn]{article}
\usepackage[a4paper, margin=1.8cm]{geometry}

\usepackage{graphicx}          % Include this line if your 
\usepackage{booktabs,float,multirow,mathrsfs,textcomp,threeparttable,siunitx,amsfonts,enumerate,algorithm,algpseudocode,xcolor}
\usepackage{amsmath,amssymb}
\usepackage{tabularx,colortbl}
\usepackage{comment, ulem}
\usepackage{mathtools}
\usepackage{hyperref}
\usepackage{tikz,pgfplots}
\definecolor{mycolor1}{rgb}{0.00000,0.44700,0.74100}%
\definecolor{mycolor2}{rgb}{0.85000,0.32500,0.09800}%
\definecolor{mycolor3}{rgb}{0.92900,0.69400,0.12500}%
\definecolor{mycolor4}{rgb}{0.49400,0.18400,0.55600}%

\newcommand{\rank}[1]{\text{rank}\left(#1\right)}
\newcommand{\norm}[1]{\left\lVert #1 \right\rVert}
\newcommand{\real}{\mathbb{R}}

\newcommand{\rto}{\rightarrow}
\newcommand{\sto}{ \rm{s.t.} }
\newcommand{\R}{\mathbb{R}}
\newcommand{\argmin}[1]{\underset{#1}{\mathrm{argmin\,}}}

\newtheorem{theorem}{Theorem}
\newtheorem{proposition}{Proposition}

\newtheorem{remark}{Remark}

\title{A constrained optimization approach to nonlinear system identification through simulation error minimization} % Title, preferably not more 
\author{Vito Cerone, Sophie M. Fosson, Simone Pirrera, Diego Regruto\thanks{The authors are with the Dipartimento di Automatica e Informatica, Politecnico di Torino, corso Duca degli Abruzzi 24, 10129 Torino, Italy; e-mail: \{vito.cerone; sophie.fosson; simone.pirrera; diego.regruto\}@polito.it}}%\\Corresponding author S. Pirrera. Tel. +39-3458029965, email: simone.pirrera@polito.it}}
\date{}

\begin{document}
\maketitle

\begin{abstract}
   This paper introduces a novel approach to system identification for nonlinear input-output models that minimizes the simulation error and frames the problem as a constrained optimization task. The proposed method addresses vanishing gradient issues, enabling faster convergence than traditional gradient-based techniques. We present an algorithm based on feedback linearization control of Lagrange multipliers and conduct a theoretical analysis of its performance. We prove that the algorithm converges to a local minimum, and it enhances computational efficiency by exploiting the problem's structure. Numerical experiments demonstrate that our approach outperforms gradient-based methods in both computational effort and estimation accuracy.
\end{abstract}

\allowdisplaybreaks

\section{Introduction}
System identification (SI) is the science of learning models of dynamical systems from experimental data. From a mathematical perspective, this process involves selecting a model class and estimating its parameters by solving an optimization problem. For further reading, see the works by \cite{ljung99}, \cite{bohlin2006practical}, and \cite{mila13}.

Most approaches to SI rely on prediction-error minimization, which involves identifying the model parameters by minimizing the difference between the measured output and the one-step-ahead prediction generated by the model.
Prediction-error minimization methods are consistent, meaning that the parameter estimates converge to the true value as the amount of data increases, under the crucial assumption that the selected model classes for both the system and the noise are accurate; see~\cite{ljung99}. 
However, if these assumptions do not hold, prediction-error minimization methods may yield models with limited accuracy in system simulation; see~\cite{zha04} and~\cite{piroddi2008simulation} for a detailed discussion of this topic.\\
A viable alternative formulates the parameter estimation problem as minimizing the simulation error, defined as the difference between the measured output and the output predicted by the model across the entire data set. This approach, known as \textit{simulation-error minimization} (SEM), yields consistent estimates regardless of the measurement noise model; see \cite{soderstrom1982some} and \cite{farina2010convergence}.

%{\color{blue} The open problem}
The primary disadvantage of SEM lies in the difficulty of solving the optimization problem, which is generically nonconvex. Notably, even when the system is linear time-invariant, the resulting SEM has a polynomial structure, thus remaining nonconvex; see, for example, \cite{cer12} and \cite{ecc23}.\\
In the general case of nonlinear systems, SEM is a differentiable, nonconvex minimization problem, and its standard solution relies on gradient-based optimization. 
Notable examples include gradient descent, pseudo-second-order algorithms, and methods based on extended Kalman filters; for a comprehensive overview, refer to \cite{ruder2016overview}. All gradient-based methods require computing the gradients of the output samples with respect to the model parameters. These gradients are recursively dependent, resulting in a dynamic relationship. \cite{Narendra91} explicitly leverage this relation to define an estimation algorithm. 
Alternatively, we can compute the gradients in closed form or through automatic differentiation; see \cite{baydin2018automatic} for a recent survey.  
Among gradient-based methods based on automatic differentiation, we mention the backpropagation through time (BPTT) algorithm, as introduced by \cite{williams1990gradient}, which has gained renewed popularity for identifying nonlinear systems modeled using recurrent neural networks. For further reading, refer to \cite{medsker1999recurrent} as well as \cite{salem2022recurrent}.

A drawback of gradient-based methods is that the dynamic relationships among gradient samples can result in vanishing and exploding gradient issues. In essence, these problem arise from the stability of the dynamical system that governs the interactions among gradients. If the system is stable, the contribution of past samples vanishes; if it is unstable, the gradients diverge. The practical implications include convergence to inaccurate estimates and either slow convergence or instability of the identification algorithm. For further details, refer to \cite{bengio1994learning} and \cite{pascanu2013difficulty}. \\
When dealing with gray-box SI, gradient-based SEM often results in highly inaccurate estimates; see, e.g., \cite{parrilo2003initialization} and \cite{mer19}. In black-box SI, the standard approach to mitigating the vanishing gradient problem involves constructing ad hoc models like long-short-term memory (LSTM) or gated recurrent unit (GRU) networks. However, LSTM and GRU tend to overfit, and their training is computationally intensive and time-consuming. 

Alternatives to gradient-based optimization are scarce. While meta-heuristic approaches, such as those proposed by \cite{blanco2001real} and \cite{bas2022training}, are available, they are often disregarded due to high computational demands and the absence of theoretical guarantees. Furthermore, studies by \cite{piroddi2003identification} and \cite{farina2011simulation} explore prediction-error minimization and a posteriori simulation-error evaluation. However, these methods exclusively apply to linearly parameterized models.
\cite{adeoye2024inexact} introduce a method based on constrained optimization and sequential quadratic programming (SQP). However, since SQP is a second-order algorithm that requires computing large Hessian matrices, approximate formulations are adopted. 

This paper proposes a novel algorithm to address SEM identification of nonlinear input-output (NIO) models. We formulate the identification problem as a \textit{constrained optimization} task. The cost function accounts for the simulation error, while the constraints enforce the model's structure. This approach circumvents recursively computing the gradients of output variables over time, thus eliminating the vanishing or exploding gradient issues, although at the cost of introducing additional optimization variables.\\
To solve the proposed optimization problem, we leverage the \textit{feedback linearization controlled multiplier optimization} (FL-CMO), a control-based continuous-time method for equality-constrained optimization first introduced by \cite{cerone2024new} and later extended by \cite{cdc24} and \cite{centorrino2025} to handle inequality constraints and non-smooth problems, respectively. 
More precisely, we derive and analyze a stable iterative algorithm from the FL-CMO differential equation, by using Euler integration.

This work introduces two primary contributions. First, it conducts a comprehensive theoretical analysis of the methodological approach. Specifically, the analysis concerns the relationship between the constrained and unconstrained SEM problem formulations and establishes convergence guarantees for the proposed method. We also examine the algorithm's computational complexity and propose a strategy to mitigate its bottlenecks. Moreover, we demonstrate how to adapt the approach to handle noisy input data and state-space models.

The second contribution demonstrates the effectiveness of the proposed approach in solving several significant SI problems. The study considers four key examples. The first example shows that the proposed method outperforms standard methods for neural NIO identification on the fluid damper benchmark. In the second example, we compare our approach with state-of-the-art results on the Bouc-Wen benchmark. The third example illustrates that the proposed method enhances accuracy and reduces training time compared to state-of-the-art neural identification techniques, such as  LSTM and GRU. Finally, in the fourth example, we address a realistic gray-box identification problem and compare our method with prediction-error minimization and gradient-based SEM.

The structure of the paper is as follows. Section \ref{sec:probl_formul} formalizes the SEM identification problem and examines the role of the vanishing gradient issue in the existing literature. Section \ref{sec:training} expands the proposed approach, formulating the problem as constrained optimization and developing the algorithm based on FL-CMO. Section \ref{sec:th_anali} provides a theoretical analysis of the developed algorithm, establishing its relation to the standard formulation and proving its convergence. In Section \ref{sec:complexity}, we leverage the problem structure and sparse QR factorization to reduce the computational complexity. Section \ref{sec:variants} extends the proposed approach to handle errors-in-variables noise and state-space models. Finally, Section \ref{sec:exampl} presents numerical examples illustrating the effectiveness of the proposed approach, and Section \ref{sec:concl} summarizes the findings and concludes the paper.

\section{Proposed System Identification Framework}
\label{sec:probl_formul}
\subsection{Problem Setup and Modeling}
Consider an unknown dynamical system of the form
\begin{equation}\label{dt_ss}\begin{aligned}
    x_{t+1} &= f(x_t,u_t) \\
    y_t &= h(x_t,u_t)
\end{aligned}\end{equation}
where $f: \R^n \times \R^q \rto \R^n$, $h: \R^n \times \R^q \rto \R^p$ are unknown functions, $n$ is the system order, $x_t \in \R^n$ is the state vector, and $u_t \in \R^q, y_t \in \R^p$ are the system's input and output of dimensions $q\in \mathbb{N}$ and $p \in \mathbb{N}$, respectively. 
As discussed, for example, by \cite{conte1999nonlinear}, observability of system \eqref{dt_ss} implies the existence of the input-output description 
\begin{equation}
    y_t = \phi(y_{t-1},\dots,y_{t-n},u_t,\dots,u_{t-n})
\end{equation}
where $\phi: \underbrace{\R^{p} \times \dots \R^p}_{n\,\rm{times}} \times \underbrace{\R^q \times \dots \R^q}_{n+1\,\rm{times}} \rightarrow \R^p$ is a suitable function depending on $f$ and $h$.

This work addresses the identification of NIO models of the form
\begin{equation}\label{eq:nio_model}
    y_t = \mathcal{M}(y_{t-1},\dots,y_{t-n},u_t,\dots,u_{t-n},\theta),
\end{equation}
where $\mathcal{M}: \underbrace{\R^{p} \times \dots \R^p}_{n\,\rm{times}} \times \underbrace{\R^q \times \dots \R^q}_{n+1\,\rm{times}} \times \R^{n_\theta} \rightarrow \R^p$ is a gray-box or black-box parametric model for $\phi$ and $\theta \in \R^{n_\theta}$ is a parameter vector of dimension $n_\theta \in \mathbb{N}$.

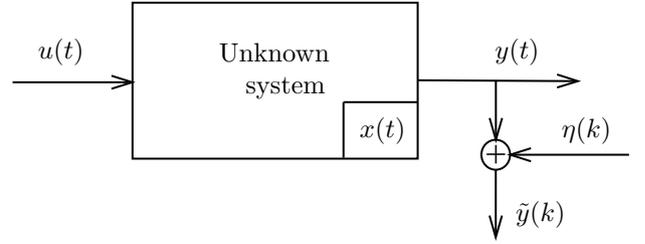
\begin{figure}
    \centering
    \tikzset{every picture/.style={line width=0.75pt}} %set default line width to 0.75pt        

\begin{tikzpicture}[x=0.75pt,y=0.75pt,yscale=-1,xscale=0.95]
%uncomment if require: \path (0,223); %set diagram left start at 0, and has height of 223

%Shape: Rectangle [id:dp24504571932551444] 
\draw   (220.14,31) -- (370,31) -- (370,109.43) -- (220.14,109.43) -- cycle ;
%Arrow [id:da49325097083382197] 
\draw    (157.14,71.01) -- (218.01,71.01) ; % detected numbers: 157.14,71.01,218.01,71.01. angle=0.0
\draw [shift={(220,71)}, rotate = 180] [color={rgb, 255:red, 0; green, 0; blue, 0 }  ][line width=0.75]    (10.93,-3.29) -- (0,0) -- (10.93,3.29)   ;
%Arrow [id:da07647861870396988] 
\draw    (370,70) -- (452.14,70.42) ; % detected numbers: 370.0,70.0,452.14,70.42. angle=0.005113176768357548
\draw [shift={(454.14,70.43)}, rotate = 180.29] [color={rgb, 255:red, 0; green, 0; blue, 0 }  ][line width=0.75]    (10.93,-3.29) -- (0,0) -- (10.93,3.29)   ;
%Straight Lines [id:da22628561604554487] 
\draw    (331,109) -- (331,84.43) -- (331,81) ;
%Straight Lines [id:da22883594824301623] 
\draw    (370.14,81) -- (331,81) ;
%Arrow [id:da9825532403515267] 
\draw    (411.02,115) -- (411,147) ; % detected numbers: 411.02,115.0,411.0,147.0. angle=1.5701713268762774
\draw [shift={(411,149)}, rotate = 270] [color={rgb, 255:red, 0; green, 0; blue, 0 }  ][line width=0.75]    (10.93,-3.29) -- (0,0) -- (10.93,3.29)   ;
%Arrow [id:da11435600018096315] 
\draw    (411.07,70.21) -- (411.07,97.93) ; % detected numbers: 411.07,70.21,411.07,97.93. angle=1.5707963267948966
\draw [shift={(411.07,99.93)}, rotate = 270] [color={rgb, 255:red, 0; green, 0; blue, 0 }  ][line width=0.75]    (10.93,-3.29) -- (0,0) -- (10.93,3.29)   ;
%Arrow [id:da2967850480448262] 
\draw    (481.04,107.46) -- (420.54,107.46) ; % detected numbers: 481.04,107.46,420.54,107.46. angle=0.0
\draw [shift={(418.54,107.46)}, rotate = 360] [color={rgb, 255:red, 0; green, 0; blue, 0 }  ][line width=0.75]    (10.93,-3.29) -- (0,0) -- (10.93,3.29)   ;
%Shape: Circle [id:dp5051590867124396] 
\draw   (403.46,107.46) .. controls (403.46,103.3) and (406.84,99.93) .. (411,99.93) .. controls (415.16,99.93) and (418.54,103.3) .. (418.54,107.46) .. controls (418.54,111.63) and (415.16,115) .. (411,115) .. controls (406.84,115) and (403.46,111.63) .. (403.46,107.46) -- cycle ;

% Text Node
\draw (257,48) node [anchor=north west][inner sep=0.75pt]    {$ \begin{array}{l}
\text{Unknown}\ \\
\ \ \ \text{system}
\end{array}$};
% Text Node
\draw (409,47.4) node [anchor=north west][inner sep=0.75pt]    {$y( t)$};
% Text Node
\draw (338,87.4) node [anchor=north west][inner sep=0.75pt]    {$x( t)$};
% Text Node
\draw (420,129.4) node [anchor=north west][inner sep=0.75pt]    {$\tilde{y}( k)$};
% Text Node
\draw (169,47.4) node [anchor=north west][inner sep=0.75pt]    {$u( t)$};
% Text Node
\draw (403.8,100.8) node [anchor=north west][inner sep=0.75pt]  [font=\small] {\textbf{+}}; %[font=\tiny]  {$\Sigma $};
% Text Node
\draw (444,87.4) node [anchor=north west][inner sep=0.75pt]    {$\eta ( k)$};

\end{tikzpicture}
    \caption{Block diagram of the data generation process}
    \label{fig:exp_setup}
\end{figure}

We assume that input-output data are available. Figure \ref{fig:exp_setup} illustrates the structure of the data-generating process. Specifically, 
The input $u_t$ is known exactly, whereas the output measurements $\tilde y_t$ are corrupted by noise according to the equation $\tilde{y}_t = y_t + \eta_t$. Here, $y_t$ represents the noise-free output, and $\eta_t \in \mathbb{R}^p$ denotes the noise sequence.

The NIO identification problem is formulated using the SEM approach; see, for example, \cite{farina2010convergence}. Accordingly, the following optimization problem is considered:
\begin{equation}\label{eq:uncons_formulation_wy0}
    \argmin{\zeta \in \R^{n_\theta+pn}} \sum_{t=1}^N \norm{ y_t(\zeta)-\tilde y_t}_W^2 + \rho(\theta)
\end{equation}
where $\zeta = [\theta^\top ,y_1^\top,\dots,y_n^\top]^\top \in \R^{n_\theta+pn}$ is the vector of optimization variables, $\norm{\eta}_W \doteq \sqrt{\eta^\top W \eta}$ for some positive definite matrix $W \in \R^{p\times p}$, $\rho: \R^{n_\theta} \rightarrow \R$ is a differentiable regularization term, and $y_t(\zeta)$ is recursively defined by 
\begin{equation}
    y_{t}(\zeta) = \mathcal{M}(y_{t-1}(\zeta),\dots,y_{t-n}(\zeta),u_t,\dots,u_{t-n},\theta), \label{def_yk_rec} 
\end{equation}
for $t=n+1,\dots,N$.

Problem \eqref{eq:uncons_formulation_wy0} generalizes the standard machine-learning formulation, as shown for example in \cite{salem2022recurrent}:
\begin{equation}\label{eq:uncons_formulation}
    \argmin{\theta \in \R^{n_\theta}} \sum_{t=1}^N \norm{y_t-\tilde y_t(\theta)}_W^2 + \rho(\theta)
\end{equation}
where $y_t(\theta)$ is the model output given known initial conditions. Typical assumptions include $y_t = \tilde y_t$ or $y_t = 0$ for $t = 1,\dots,n$ in input-output models, and $x_0=0$ or a random normally distributed value for state-space models. Adopting formulation \eqref{eq:uncons_formulation_wy0} enables joint estimation of the initial conditions and the network parameters, leading to more accurate estimates.

\subsection{Background and Related Work}
\label{sec:related_literature}
The standard approach to solving Problems \eqref{eq:uncons_formulation_wy0} and \eqref{eq:uncons_formulation} is to apply gradient-based methods. Several techniques are available in this context. The most popular approach involves gradient-descent (GD) algorithms, such as stochastic and mini-batch GD, RMSprop, and Adam; see \cite{ruder2016overview} for a comprehensive review. Alternatively, second-order and pseudo-second-order methods are employed. Popular choices include Newton's method, Broyden-Fletcher-Goldfarb-Shanno algorithms (see \cite{peng2011nonmonotone}, \cite{liu2018limited}, and \cite{bemporad2025bfgs}), and methods tailored to nonlinear least-squares problems such as the Gauss-Newton or Levenberg–Marquardt algorithms (see \cite{mirikitani2009recursive}). 

All the approaches mentioned above require computing the gradients $\nabla_\theta y_t$, which recursively depend on $\nabla_\theta y_{j}$, for all $j=1,\dots,t-1$. This dependence implicitly defines a dynamic relation between the gradient samples. If the dynamics is stable, the contribution to $\nabla_\theta y_t$ due to $\nabla_\theta y_k$, for $k \ll t$, is negligible (vanishing gradient). Conversely, $\nabla_\theta y_t$ is divergent if the dynamics is unstable. These phenomena are known as the vanishing and exploding gradient issues, respectively. These, discussed in \cite{bengio1994learning,pascanu2013difficulty}, cause difficulties in the learning process such as slow convergence, inaccuracy, and numerical instability.

A practical method to mitigate the exploding-gradient issue is gradient clipping, as discussed by \cite{ramaswamy2023gradient} and the references therein. However, the vanishing-gradient problem remains unresolved in general. In black-box SI, a standard solution is to use LSTM or GRU networks, as noted by \cite{salem2022recurrent}. These networks are designed to reduce the impact of the vanishing-gradient problem, although they do not fully resolve it. Alternatively, existing gradient-free algorithms, such as that proposed in~\cite{bas2022training}, are unaffected by vanishing and exploding gradients but converge slowly and require significant computational resources.
\cite{adeoye2024inexact} suggest formulating the identification problem as a constrained optimization task to develop a learning algorithm for recurrent neural networks using SQP. They argue that this formulation enables the development of alternative training algorithms for neural networks and provides tools to derive convergence results. 

In this work, we explore a constrained optimization approach, demonstrating that it not only leads to the development of new algorithms but also ensures their resilience to vanishing-gradient issues.

\section{Identification using FL-CMO algorithm}
\label{sec:training}
This section introduces a novel algorithm for the identification of NIO models. The method is based on the FL-CMO algorithm developed in \cite{cerone2024new}. 

\subsection{Constrained optimization {formulation}}
\label{sec:cns_opt_appr}

In this section, we reformulate Problem \eqref{eq:uncons_formulation_wy0} as a \textit{constrained optimization} problem, i.e., we solve:
\begin{align}\label{cns_opt_nio}
    &\argmin{y_1 \in \R^{p},\dots,y_{N} \in \R^p, \theta \in \R^{n_\theta}} \sum_{t=1}^N \norm{ \tilde{y}_t - y_t }_W^2  + \rho(\theta) \notag \\
    & \quad \textrm{s.t.} ~~ y_t = \mathcal{M}(y_{t-1},\dots,y_{t-n}, u_t,\dots, u_{t-n},\theta)  \\
    & \qquad \text{for } t = n+1,\dots,N.\notag 
\end{align}
Similar mathematical formulations can be found in \cite{adeoye2024inexact}, who introduce a method based on inexact SQP for training RNNs. Analogous constrained optimization problems are also formulated in the set-membership identification setting (see, e.g., \cite{milanese2004set}, \cite{mila13}, \cite{cer14}), under the typical assumption that the noise is bounded and incorporates additional inequality constraints accordingly.

We denote $\nu \doteq {n_\theta + pN}$, $m\doteq N-n$, and
\begin{subequations}\begin{align}
    \xi &= [\theta^\top, y_1^\top, \dots, y_N^\top]^\top \in \R^\nu,\\
    h_t(\xi) &= y_{t+n}-\mathcal{M}(y_{t+n-1},\dots,y_{t}, u_{t+n},\dots, \notag \\
    & \quad u_{t},\theta), ~~t=1,\dots,m,\\
    \mathcal{L}(\xi) &= \sum_{t=1}^N \norm{ \tilde{y}_t - y_t }_W^2  + \rho(\theta)
\end{align}\end{subequations}
the optimization variables, constraints, and cost function of Problem \eqref{cns_opt_nio}, respectively. Let $\lambda \in \R^m$ be the vector of Lagrange multipliers. The first-order optimality conditions for Problem \eqref{cns_opt_nio} are:
\begin{subequations}\label{eq:kkt_cns}
\begin{align}
	& \nabla_{\xi} \mathcal{L}(\xi) + \sum_{t=1}^m \lambda_{t} \nabla_\xi h_t(\xi) = 0\\
	& h_{t}(\xi)=0, \qquad t = 1,\dots,m.
\end{align}
\end{subequations}
\begin{remark}
    Conditions \eqref{eq:kkt_cns} require computing only the gradients of simple functions that do not recursively depend on other expressions. As a result, solving \eqref{eq:kkt_cns} avoids the numerical issues related to the vanishing and exploding gradient problems.
\end{remark}

\subsection{A motivating example}
\label{sec:motive_example_part2}
  To substantiate the necessity of the proposed method, we present an example illustrating how the vanishing and exploding gradient issues occur in a basic scenario and how one can address them through  constrained optimization.

Consider the following first-order linear-time-invariant model:
\begin{equation}\label{eq:ex_system}y_t = \theta_1 y_{t-1} + \theta_2 u_{t-1}, \qquad y_1=0.
\end{equation}
The model output sample at time $t$ is given by:
%$$y_{t} = \alpha^{t-1} y_1 + \sum_{k=1}^{t-1} \alpha^{k-1} \beta u_{t-k}. $$
\begin{equation} y_{t} = \sum_{k=1}^{t-1} \theta_1^{k-1} \theta_2 u_{t-k}.
\end{equation}
Problem \eqref{eq:uncons_formulation} with $W=1$ and $\rho(\theta)=0$ becomes
%$$\mathcal{L}(\alpha,\beta,y_1) = \sum_{t=2}^N \left( \tilde y_t - \alpha^{t-1} y_1 - \sum_{k=1}^{t-1} \alpha^{k-1} \beta u_{t-k} \right)^2 $$
\begin{equation}{\min_{\theta \in \R^{n_\theta}}}\,\mathcal{L}(\theta) = 
%\sum_{t=2}^N \left( y_t - \tilde y_t\right)^2 =
\min_{\theta \in \R^{n_\theta}} \sum_{t=2}^N \left( -\tilde y_t + \sum_{k=1}^{t-1} \theta_1^{k-1} \theta_2 u_{t-k} \right)^2, \end{equation}
and its optimal solution satisfies the optimality condition:
\begin{equation}\label{eq:motiv_ex_part1a}\nabla_{\theta} \mathcal{L} = 2 \sum_{t=2}^N \left(  y_t -\tilde y_t\right) \nabla_\theta y_t {=0,} \end{equation}
where
\begin{equation} \label{eq:motiv_ex_part1b}
\nabla_\theta y_t = \begin{bmatrix} \frac{d y_t}{d \theta_1} \\ \frac{d y_t}{d \theta_2} \end{bmatrix} = \sum_{k=1}^{t-1} \begin{bmatrix} (k-1)\theta_1^{k-2}\theta_2 u_{t-k} \\ \theta_1^{k-1} u_{t-k} \end{bmatrix}. 
\end{equation}
In a typical gradient-descent iteration, if $ \vert \theta_1 \vert < 1 $, the terms in the summation associated with large values of $k$ become negligible, which can lead to numerical cancellation issues. As a result, for each time step $t$, only a few terms in the sum from Equation \eqref{eq:motiv_ex_part1b} significantly impact the update of the parameter estimate. On the other hand, if $ \vert \theta_1 \vert > 1 $, the norm $\norm{ \nabla_\theta y_t } $ grows exponentially as $t$ increases.

Using formulation \eqref{cns_opt_nio} and according to the conditions in \eqref{eq:kkt_cns}, the optimality conditions of the problem are:
\begin{subequations}\label{eq:motiv_ex_part2}
\begin{equation}
   \sum_{t=2}^N \lambda_{t-1} y_{t-1}=0, \qquad \sum_{t=2}^N \lambda_{t-1} u_{t-1}=0 
 \end{equation}
\begin{equation}
 2 (y_t - \tilde y_t) +\lambda_{t-1} - \lambda_{t} \theta_{1} = 0, \quad \forall t = 1,\dots,N, 
\end{equation}
\begin{equation} 
h_{t-n}= y_t-\theta_1 y_{t-1}-\theta_2 u_t=0, \quad \forall t = 2,\dots,N,
\end{equation}\end{subequations}
where $\lambda_{N}=0$. 

Compared to \eqref{eq:motiv_ex_part1a}, the system of equations \eqref{eq:motiv_ex_part2} 
%contains more variables and equations:
has larger dimensions: the number of variables increases from $2$ to $2+N$, and the number of equations from $2$ to $2N+1$. However, all equations are now simple bilinear functions of the unknowns $x$ and $\lambda$, eliminating the computation of the  powers of $\theta_1$, and thus circumventing vanishing and exploding gradient problems.

%In general, the formulation \eqref{cns_opt_nio} avoids vanishing and exploding-gradient phenomena. The novel structure of the problem does not involve iteratively composing of the function $\mathcal{M}$, which previously led to repeated multiplication of the gradients by the Jacobian of $\mathcal{M}$ when using gradient-based methods to solve Problem \eqref{eq:uncons_formulation_wy0}.

%\subsection{Review of FL-CMO algorithm}
\subsection{Identification algorithm based on FL-CMO}
This section briefly reviews the FL-CMO dynamics introduced by \cite{cerone2024new} and presents an iterative algorithm based on its Euler discretization.

FL-CMO addresses equality-constrained optimization problems of the form
\begin{equation}\label{eq_constr}\begin{aligned}
    \min_{\xi \in \R^\nu}& ~f(\xi) \\
    &\textrm{s.t.} \quad h(\xi)=0,
\end{aligned} \end{equation}
where $f: \mathbb{R}^\nu \rightarrow \mathbb{R}, h: \mathbb{R}^\nu \rightarrow \mathbb{R}^m$ are differentiable, possibly non-convex functions. The dynamics is developed using design techniques in the context of continuous-time feedback control systems, and the following equations describe the plant to be controlled:
\begin{subequations}\label{plant_cmo}
\begin{align}
\dot{\xi} &= -\nabla f(\xi) - J_h(\xi) \lambda \\
y &= h(\xi).
\end{align}
\end{subequations}
The Lagrange multipliers $\lambda$ act as the input to the plant, while the optimization variables $\xi$ represent the state, and the constraints define the outputs. This definition is grounded in the observation that an equilibrium point $(\xi^*, \lambda^*)$ of \eqref{plant_cmo} is a stationary point of \eqref{eq_constr} if and only if $h(\xi^*) = 0$; see \cite[Lemma 1]{cerone2024new}. Therefore, any feedback control system that ensures asymptotic stability and zero output regulation will inherently lead to convergence to a stationary point of \eqref{eq_constr}.

Under the assumption that the Jacobian of the constraints $J_h(\xi)$ is full rank for all $\xi \in \real^\nu$, \eqref{plant_cmo} has vector relative degree $r = [1,\dots,1]^\top \in \real^m$ \cite[Lemma 2]{cerone2024new}, and therefore the non-interacting control problem defined in the context of feedback linearization admits a solution, as shown in \cite{Isidori1995} and \cite{kha02}. Applying the technique of output-feedback linearization controller design, we obtain the controller:
\begin{subequations}\label{cmo_controller}\begin{align}
    \lambda &= (J_h(\xi) J_h^\top(\xi))^{-1}(J_h(\xi) \nabla_\xi f(\xi) + v)\\
    v &= K y,  K > 0.
\end{align}\end{subequations}
The continuous-time closed-loop system \eqref{plant_cmo}-\eqref{cmo_controller} defines the FL-CMO {dynamics}, which locally converges to a solution of the constrained optimization problem \eqref{eq_constr}. We refer the reader to Section \ref{sec:converg} for a detailed convergence analysis.

We integrate the differential equations describing the FL-CMO dynamics using the Euler method; see \cite{wanner1996solving} and \cite{Sundnes2024}.
Algorithm \ref{alg:cap} 
%"below" SOP: non mettere mai sopra/sotto.. non serve! e poi non sappiamo mai dove latex mette le cose 
illustrates the resulting procedure.

\begin{algorithm}
%\caption{{Identification algorithm based on FL-CMO dynamics and Euler discretization method}}\label{alg:cap} %%SOP: caption inutilmente verbosa
\caption{{Identification algorithm based on FL-CMO  and Euler discretization}}\label{alg:cap}
\begin{algorithmic}
\Require $\xi_0,\epsilon_f,\epsilon_h,K,\tau$
\State {$k \gets 0$}
\While {$\norm{\delta_\xi}_2 \geq \epsilon_f$ \text{or} $\norm{h(\xi_k)}_2 \geq \epsilon_h$} 
\State $J = J_h(\xi_k)$
\State $\delta_\xi = -\nabla f(\xi_k) - J^\top (J J^\top)^{-1} (J \nabla f(\xi_k) +  K  h(\xi_k) )$
\State $\xi_{k+1} = \xi_{k}  + \tau \delta_\xi$
\State $k \gets k+1$
\EndWhile
\end{algorithmic}
\end{algorithm}

The computation of the Jacobian $J = J_h(\xi_k)$ of the constraints $h$ can be performed using automatic differentiation, see, e.g., \cite{baydin2018automatic}, or closed-form expressions derived through tensor calculus rules, see, e.g., \cite{laue2018computing}. Both approaches efficiently perform this computation by leveraging parallelization.

\section{Theoretical Analysis}
\label{sec:th_anali}
This section analyzes the theoretical properties of Algorithm~\ref{alg:cap}.

\subsection{Relation with the unconstrained formulation}
By construction, the constrained formulation \eqref{cns_opt_nio} and the unconstrained one \eqref{eq:uncons_formulation_wy0} share the same global optimal solutions. However, local solutions may differ between the two. The following theorem establishes a connection between the local solutions of \eqref{cns_opt_nio} and those of \eqref{eq:uncons_formulation_wy0}.

\begin{theorem}
\label{th:rel_stat}
    If a point $[\theta^\top, y_1^\top, \dots, y_N^\top]^\top \in \R^{n_\theta+pN}$ satisfies the first-order optimality conditions of Problem \eqref{cns_opt_nio}, the point $[\theta^\top, y_1^\top, \dots, y_n^\top]^\top \in \R^{n_\theta+pn} $ is a stationary point of Problem \eqref{eq:uncons_formulation_wy0}.
\end{theorem}

\textit{Proof.}
    We denote $z = [\theta^\top, y_1^\top,\dots,y_{n}^\top]^\top \in \R^{n_{\theta}+pn}$, and $w = [y_{n+1}^\top,\dots,y_{N}^\top]^\top \in \R^m$, where $m=p(N-n)$.
     Let $h(z,w)$ denote the constraints of \eqref{cns_opt_nio} and $J_h(z,w) = [J_{h,z}, J_{h,w}]$ its Jacobian. By construction, $J_{h,w} \in \R^{m\times m}$ is block-triangular with all the blocks on the main diagonal being full rank, making it invertible. Denote $G(w) = \sum_{t=n+1}^N \norm{y_t - \tilde y_t}_{W}^2$ and $R(z) = \sum_{t=1}^n \norm{y_t - \tilde y_t}_{W}^2 + \rho(\theta)$. The first-order optimality conditions are given by
    \begin{subequations}\begin{align}
        \nabla_z R + J_{h,z}^\top \lambda = 0 \label{eq:pr1_opt1}\\
        \nabla_w G + J_{h,w}^\top \lambda = 0.\label{eq:pr1_opt2}
    \end{align}\end{subequations}
    From \eqref{eq:pr1_opt2} we obtain
    \begin{equation}
        \lambda = - \left(J_{h,w}^{\top}\right)^{-1} \nabla_w G \label{eq:pr1_lamb}
    \end{equation}
    and substituting \eqref{eq:pr1_lamb} into \eqref{eq:pr1_opt1} we obtain 
    \begin{equation}\label{eq15}
        \nabla_z R - J_{h,z}^\top \left(J_{h,w}^{\top}\right)^{-1} \nabla_w G =0.
    \end{equation}
    By applying the Implicit Function Theorem, and noting that $h(z,w)=0$ and $J_{h,w}$ is invertible, there exists a function $\gamma: \R^{n_{\theta}+pn} \rightarrow \R^m$ such that $w = \gamma(z)$. Moreover, the Jacobian $J_{\gamma,z}$ of $\gamma$ satisfies
    \begin{equation}\label{jacob_imp_fun}
        J_{\gamma,z} = -J_{h,w}^{-1} J_{h,z}.
    \end{equation}
    We now consider Problem \eqref{eq:uncons_formulation_wy0}, for which the stationary points satisfy
    \begin{equation}
        \nabla_z G(\gamma(z)) + \nabla_z R(z)=0.
    \end{equation}
    By applying the chain rule and Equation \eqref{jacob_imp_fun}, we derive the following expression:
    \begin{equation}\begin{aligned}
        &\nabla_z R(z) + \nabla_z G(\gamma(z))= \\
        &= \nabla_z R + \left[\nabla_{\gamma(z)} G^\top J_{\gamma,z} \right]^\top = \\
        &= \nabla_z R -\left[ \nabla_w G^\top J_{h,w}^{-1} J_{h,z}\right]^\top = \\
        &=  \nabla_z R - J_{h,z}^\top \left(J_{h,w}^{\top}\right)^{-1}  \nabla_w G,
    \end{aligned}    \end{equation}
    which equals zero by \eqref{eq15}. $\quad \square$

Theorem \ref{th:rel_stat} emphasizes that solving \eqref{cns_opt_nio} instead of \eqref{eq:uncons_formulation_wy0} does not introduce any additional local solutions. 

\subsection{Analysis of local convergence}
\label{sec:converg}
In the following theorem, we analyze the convergence of Algorithm \ref{alg:cap}. 

\begin{theorem}[Convergence of Algorithm \ref{alg:cap}]
    Any point $(\xi^*,\lambda^*)$ satisfying the second-order sufficient conditions of Problem \eqref{cns_opt_nio} is a locally asymptotically stable equilibrium of Algorithm \ref{alg:cap} if $\tau$ is sufficiently small.   
\end{theorem}

\textit{Proof.}
    First, we note that the assumptions required to apply the FL-CMO method hold:
    \begin{itemize}
        \item By construction, Problem~\eqref{cns_opt_nio} is characterized by $m < \nu$, meaning that the number of optimization variables is larger than the number of constraints.
        \item $\rank{J_h(\xi)}=m$ for any $\xi \in \R^\nu$. This condition holds because $J_h(\xi)$ consistently maintains full rank, as demonstrated by the argument used for the $J_{h,w}$ submatrix in the proof of Theorem \ref{th:rel_stat}.
    \end{itemize}
    
    Under these conditions, Theorem 4 in \cite{cerone2024new} is applicable. Consequently, any point $(\xi^*,\lambda^*)$ satisfying the second-order sufficient conditions of Problem \eqref{cns_opt_nio} is a locally asymptotically stable equilibrium point of the dynamics \eqref{plant_cmo}-\eqref{cmo_controller}. Finally, the argument remains valid after the Euler discretization because it preserves the equilibrium points and their stability properties if $\tau$ is sufficiently small; see, e.g.,~\cite{wanner1996solving}.  $\quad \square$

{The assumption that $\tau$ must be sufficiently small to ensure stability closely parallels the requirement for a small learning rate in standard gradient-based optimization for unconstrained problems in machine learning. As in those settings, the value of $\tau$ must be chosen carefully: if it is too large, the algorithm diverges; if it is too small, convergence becomes slow.}

\section{Analysis of the computational complexity}
\label{sec:complexity}
The primary computational cost of Algorithm \ref{alg:cap} arises from computing:
\begin{equation}\label{lin_sys}
    {\sigma \doteq (J_h(\xi_k) J_h^\top(\xi_k))^{-1} (J_h(\xi_k) \nabla f(\xi_k) + K h(\xi_k)).}
\end{equation}
%{\color{magenta} **** qui non è chiaro: \label{lin_sys} abbiamo un sistema di equazioni lineari o la **soluzione** di un sistema di equazioni lineari? Chi è $\sigma$? Inoltre va aggiunta una frase per spiegare a parola il core della complessità, che sta nell'inversione... ma va detto.}

This section demonstrates how the problem's structure can be leveraged to reduce the computational complexity related to the main computational bottleneck of Algorithm~\ref{alg:cap}: the inversion of the matrix $(J_h(\xi_k) J_h^\top(\xi_k))$. 

Suppose that we use a standard decomposition method to compute the triangular Cholesky factor $R$ of $J_h(\xi_k) J_h^\top(\xi_k)$. The computational complexity of such an operation is $O(m^3)$. See \cite{stewart1998matrix} for a detailed analysis of the computational complexity of these algorithms.
In the following, we present a method for solving the linear system \eqref{lin_sys} that exploits the sparsity of $J_h(\xi_k)$ to reduce the overall computational effort to $O(m^2)$. Specifically, the method proceeds as follows:
\begin{enumerate}
    \item Perform the Q-free QR decomposition of $J_h^\top(\xi_k)$, i.e., obtain the upper triangular factor $R \in \real^{m \times m}$, where an orthogonal matrix $Q$ satisfies $J_h^\top(\xi_k) = Q R$. This method employs the sparse Householder algorithm described in \cite{davis2011algorithm}. By applying the QR factorization, we obtain $J_h(\xi_k) J_h^\top(\xi_k) = R^\top Q^\top Q R = R^\top R$, which makes $R^\top$ a Cholesky factor of $J_h(\xi_k) J_h^\top(\xi_k)$. 
    \item Apply forward and backward substitution to compute the solution of the linear system.
\end{enumerate}

%In the remainder of this section, \red{we define $\phi$ as the dynamical order of the NIO model $\rightarrow$ mettere $\nu$ o $\phi$ dall'inizio.} and $n = n_\theta + pN$ as the number of optimization variables {of Problem \eqref{cns_opt_nio}.}

\begin{proposition}\label{prop:complexity}
    Using the sparse Householder algorithm, the computational complexity of a generic iteration of Algorithm \ref{alg:cap} is $O(m^2)$.
\end{proposition}

\textit{Proof}. We calculate the operations required for the Householder algorithm that exploits sparsity by avoiding multiplications involving zeros. Table \ref{tab:comp_complex} reports the number of floating-point operations (FLOPs) for each step of the Householder algorithm applied to the matrix $ X = J_h^\top(\xi_k)$, comparing the dense and sparse cases.
\begin{table*}[ht]
    \centering
    \caption{FLOPs count for each instruction of the Housholder algorithm.}
    \label{tab:comp_complex}
    \begin{tabular}{l|c|c|c}
        Instruction & Iterations & FLOPs - dense & FLOPs - sparse \\
        \hline
        $U=0_{\nu \times m}$,$R=0_{m\times m}$ & 1 & - & - \\
        for $k=1,\dots,m$ & & & \\       $\quad$ $U_{k:\nu,k},R_{k,k}=\text{housegen}(X_{k:\nu,\nu})$ & $m$ & $3(\nu-k+1)$ & $3 \chi(k)$ \\
        $\quad$ $v=U^\top_{k:\nu,k}X_{k:\nu,k+1:m}$ & $m$ & $(\nu-k+1)(m-k)$ & $\psi(k)$ \\
        $\quad$ $X_{k:\nu,k+1:m}=X_{k:\nu,k+1:m} -U_{k:\nu,k}v$ & $m$ & $(\nu-k+1)(m-k)$ & $\zeta(k)$ \\
        $\quad$ $R_{k,k+1:m} = X_{k,k+1:m}$ & $m$ & -  & -\\
        end for & & &\\
        \hline
    \end{tabular}
\end{table*}

To analyze the first instruction of the loop, we examine the function $[u, \mu] =$ \textit{housegen}$(w)$. This function performs the Householder reflection of its argument by computing $u$ such that $(I - u u^\top)w = \pm \norm{w}_2 e_1$, where $e_1$ is the first vector of the standard Euclidean basis. For additional details, refer to \cite{stewart1998matrix}.
\begin{algorithm} 
\caption{The \textit{housegen} function}
\label{alg:housegen} 
\begin{algorithmic} 
\State $u = w, \quad \mu = \norm{w}_2$
\If{$\nu = 0$}\; \State $u_1 = \sqrt2,\quad$ {return} \; \EndIf 
\State $u = w / \mu$
\If{$u_1 \geq 0$} 
\State $u_1 = u_1 + 1 ,\quad \mu = -\mu $
\Else 
\State $u_1 = u_1 - 1 $
\EndIf 
\State $u = u / \sqrt{\vert u_1\vert}$
\end{algorithmic} 
\end{algorithm}
The \textit{housegen} function is defined according to Algorithm \ref{alg:housegen}. The computational cost of \textit{housegen} is primarily driven by the calculation of $\mu=\norm{w}_2$, $u = w/\mu$, and $u = u / \sqrt{\vert u_1\vert}$. Each of these operations involves $\chi$ FLOP products, where $\chi$ is the cardinality of $w$. Consequently, since $\xi_k$ is constructed from the $k$-th to the $n$-th element of $X$, this step requires $3(\nu-k+1)$ FLOPs in the dense case and $3\chi(k)$ FLOPs when exploiting sparsity. Based on the structure of $X = J_h^\top(\xi_k)$, $\chi(k)$ is given by 
\begin{equation}
    \chi(k) = \begin{cases} n_\theta -k+1+p(1+n) &\quad \text{ if } k \leq n_\theta \\ p(1+n) &\quad\text{ if } k > n_\theta. \end{cases}
\end{equation}
Computing the product $U_{k:\nu,k}^\top X_{k:\nu,k+1:m}$ in the dense case requires $(\nu-k+1)(m-k)$ multiplications. In the sparse case, however, $U_{k:\nu,k}$ exhibits the same sparsity pattern as $\xi_k$, except for its first element, which consistently provides no contribution to the product. By leveraging this sparsity, the number of required products drops to
\begin{equation}
\psi(k) = \begin{cases}\psi_1(k)  &\quad \text{ if } k \leq n_\theta \\ 
    \psi_2(k) &\quad\text{ if } k > n_\theta,
\end{cases} 
\end{equation}
where, for $i = (k \text{ mod } p)+1$, we obtain 
\begin{equation}\label{psi2} \begin{aligned}
    \psi_2(k)& = (p-i)p(n+1)+n p^2+\dots+2p^2+p^2\\
    &=(p-i)p(n+1)+\sum_{z=1}^n z p^2,
\end{aligned}\end{equation}
\begin{equation}\label{psi1}\begin{aligned}
    \psi_1(k) &= (n_\theta-k+1)(m-k) + \psi_2(k).
\end{aligned}\end{equation}
Equation \eqref{psi2} follows from the observation that the number of coinciding nonzero indices decreases by $p$ every $p$ columns. Similarly, Equation \eqref{psi1} is obtained using the same procedure while also accounting for the contributions of products involving the first $(n_\theta-k+1)$ elements of $U_{k:\nu,k}^\top$.

We now compute $\zeta(k)$ representing the number of FLOPs required to calculate the product $U_{k:\nu,k} v$, where $U_{k:\nu,k} \in \R^{\nu-k+1}$ is a sparse vector with cardinality
\begin{equation}
    \chi_1(k) = \begin{cases}
        \chi(k) & \quad \text{ if } k\leq n_\theta \\
        \chi(k)+1 & \quad \text{ if } k > n_\theta
    \end{cases}
\end{equation}
and $v \in \R^{1,m-k}$ is possibly dense. Accordingly, the result is an $(\nu-k+1) \times (m-k)$ matrix in which all rows corresponding to zero elements of $U_{k:\nu,k}$ are zero. As a result, the dense case requires $(\nu-k+1)(m-k)$ FLOPs, whereas the sparse case requires
\begin{equation}
    \zeta(k) = (m-k)\chi_1(k) \text{ FLOPs}.
\end{equation}
To conclude the proof, we demonstrate the existence of constants $M>0$ and $m_0>0$ such that for all $m>m_0$:
\begin{equation} \label{eq:proof_complex_maincond}
    \left\lvert \sum_{k=1}^m 3\chi(k)+\psi(k)+\zeta(k) \right\rvert \leq M m^2.
\end{equation}
Let us evaluate each term of the sum:
\begin{equation}\label{sum_chi}\begin{aligned}
    &\sum_{k=1}^m 3\chi(k) = 3 \sum_{k=1}^{n_\theta} \left( n_\theta -k+1+p(1+n) \right) +\\
    & +3 \sum_{k=n_\theta+1}^{m}p(1+n) = 3m(1+n)p+\frac{3}{2}(n_\theta^2 + n_\theta)
\end{aligned} \end{equation}
\begin{equation}\label{sum_psi}\begin{aligned}
    &\sum_{k=1}^m \psi(k) = \sum_{k=1}^{n_\theta} \psi_{1}(k) + \sum_{k=n_\theta+1}^m \psi_2(k) \\
    &\quad = \sum_{k=1}^{n_\theta} (n_\theta-k+1)(m-k) + \sum_{k=1}^m \psi_2(k)
\end{aligned} \end{equation}
\begin{equation}\label{sum_zeta}\begin{aligned}
    &\sum_{k=1}^m \zeta(k) = \sum_{k=1}^{n_\theta} (n_\theta -k+1+p(n+1))(m-k) + \\
    &\quad +\sum_{k=n_\theta+1}^{m} (p(n+1)+1)(m-k) = \\
    &= \sum_{k=1}^{n_\theta} (n_\theta -k)(m-k) +\sum_{k=1}^{m} (p(n+1)+1)(m-k) = \\
    &= \frac{m}{2} (n_\theta^2 - n_\theta)- \frac{n_\theta}{6}(n_\theta^2 -1) + \frac{m^2-m}{2}(p(n+1)+1)\\
    &= \frac{m^2}{2}(p(n+1)+1) + \\
    & \quad +\frac{m}{2} (n_\theta^2 - n_\theta-p(n+1)-1)- \frac{n_\theta}{6}(n_\theta^2 -1).
\end{aligned}\end{equation}
Let us expand individually the two terms in equation \eqref{sum_psi} to analyze their contributions.
\begin{equation}\begin{aligned}
    & \sum_{k=1}^{n_\theta} (n_\theta-k+1)(m-k) = \\
    &= \frac{m}{2} (n_\theta^2 - n_\theta)- \frac{n_\theta}{6}(n_\theta^2 + 3n_\theta + 2)
\end{aligned}\end{equation}
\begin{equation}\begin{aligned}
    & \sum_{k=1}^{m} \psi_{2}(k) = \sum_{k=1}^{m}\left( (p-i)(n+1)p + \sum_{z=1}^n p^2 z \right) = \\
    &\quad = \sum_{j=1}^{m/p} \sum_{i=1}^p \left( (p-i)(n+1)p + \sum_{z=1}^n p^2 z \right) =\\
    &\quad = \frac{m}{p} \left( p^3 (n+1) - \frac{1}{2}(n+1)p^2(p+1) \right) + \\
    & \qquad+\frac{1}{2} m p^2 n (n+1) = \frac{1}{2} mp(n + 1)(p + p n - 1)
\end{aligned}\end{equation}
From the above equations, we observe that \eqref{sum_zeta} dominates \eqref{sum_chi} and \eqref{sum_psi}. By choosing $M$ in equation \eqref{eq:proof_complex_maincond} to be greater than $\frac{1}{2}+\frac{1}{2}p(n+1)$, the result follows. $\quad \square$

\begin{remark}
    The sparsity of $J_h(\xi_k)$ can also be exploited when computing the matrix-vector products $J_h(\xi_k) \nabla f(\xi_k)$ and $J_h^\top(\xi_k) \sigma$. Table \ref{tab:flops_mult} presents the FLOP count for both dense and sparse algorithms.
\end{remark}

\begin{table}[ht]
    \centering
    \caption{FLOP count for matrix-vector multiplications}
    \label{tab:flops_mult}
    \begin{tabular}{l|c|c}
        Instruction &  FLOPs - dense & FLOPs - sparse \\
        \hline
        $J_h(\xi_k) \nabla f(\xi_k)$ & $m(n_\theta+N)$ & $m(n_\theta + pn)$ \\
        $J_h^\top(\xi_k) \sigma$ & $m^2$ & $m(n_\theta + pn)$\\
        \hline
    \end{tabular}
\end{table}

We evaluate the practical implications of the achieved reduction in computational complexity proven in Proposition~\ref{prop:complexity} through a series of numerical experiments. We construct a synthetic problem with $p=n=2$ and measure the time required to perform a single iteration of Algorithm \ref{alg:cap} across varying dataset dimensions. The results demonstrate relative speedups of $1.6\times, 3.7\times$, and $5.7\times$ for $N=10^3, N=5\cdot10^3$, and $N=10^4$, respectively.

Compared to the inexact SQP method proposed by \cite{adeoye2024inexact}, we note that the size of the matrix requiring QR factorization decreases from $(\nu+m) \times (\nu+m)$ to $m \times (\nu+m)$. As a consequence, we conclude that Algorithm~\ref{alg:cap} benefits from reduced computational complexity compared to the inexact SQP method by \cite{adeoye2024inexact}. Moreover, the iterations are exact, significantly simplifying the analysis.

\section{Extensions to error-in-variables and state-space models}
\label{sec:variants}
{This section demonstrates the application of the constrained\textendash optimization reformulation and Algorithm~\ref{alg:cap} to identification problems with a structure differing from that considered in Section~\ref{sec:training}.}

{\subsection{Errors-in-variables models}}
The constrained optimization approach generalizes to scenarios in which both the input and output data are corrupted by noise. In this case, the available dataset is $\{\tilde u_t,\tilde y_t\}$, for $t=1,\dots, N$, where $\tilde u_t = u_t + \epsilon_t$, $\tilde y_t = y_t + \eta_t$. Here, $\epsilon_t\in \R^q,\eta_t\in \R^p$ denote the noise components. This setting is commonly referred to as the errors-in-variables (EIV) model.

We formalize this task by solving the following optimization problem
\begin{equation}\begin{aligned}
    & \argmin{\theta \in \R^{n_{\theta}},u_1 \in \R^q, \dots ,y_N \in \R^{p}} \sum_{t=1}^N \norm{y_t - \tilde y_t }_{W_y}^2 + \\
    & \qquad + \sum_{t=1}^N \norm{u_t-\tilde{u}_t}_{W_u}^2 + \rho(\theta) \\
    & \quad \text{s.t.} \quad y_t = \mathcal{M}(y_{t-1},\dots,y_{t-n}, u_t,\dots,u_{t-n},\theta) \\
    & \qquad \qquad \text{for } t = n+1,\dots,N.
\end{aligned}\end{equation}
where $W_y, W_u \succ 0$ are weighting matrices. When output noise is expected to exceed input noise, $W_u \succ W_y$. As $W_u$ approaches infinity ($W_u \rightarrow \infty$), $\tilde u$ becomes equivalent to $u$ ($\tilde u = u$), recovering the standard output error problem.

The theoretical analysis presented in Section \ref{sec:th_anali} also applies to this case, whereas the considerations on computational complexity are different from those in Section \ref{sec:complexity} due to the distinct sparsity of the problem. Compared to Problem \eqref{cns_opt_nio}, introduced in Section \ref{sec:training} for the output-error case, the number of optimization variables rises, while the number of constraints remains unchanged. Consequently, the dimensions of the linear system \eqref{lin_sys} remain unchanged. Furthermore, $J$ is still sparse, and numerical experiments demonstrate that the proposed QR decomposition offers computational benefits. 

{\subsection{State-space models}}
\label{sec:extension_others}
State-space models are a widely used description for dynamical systems owing to their versatility, as they encompass linear systems, gray-box models, recurrent neural networks, and neural state-space models. In this section, we propose a constrained optimization formulation of the SEM problem for state-space models.

Given a dynamical model for system~\eqref{dt_ss} defined by 
\begin{equation}\begin{aligned}
    x_{t+1} = \mathcal{M}_1(x_t, u_t, \theta),\\
    y_{t} = \mathcal{M}_2(x_t, u_t, \theta)
\end{aligned}\end{equation}
where $\mathcal{M}_1: \R^{n} \times \R^{q} \times \R^{n_\theta} \rightarrow \R^n$ and $\mathcal{M}_2: \R^{n} \times \R^{q} \times \R^{n_\theta} \rightarrow \R^p$, we can recast its SEM identification as the following constrained optimization problem:
\begin{equation}\label{opt_nnss_train}\begin{aligned}
    & \argmin{\theta \in \R^{n_\theta},x_1 \in \R^{n},\dots,x_{N} \in \R^n} \sum_{t=1}^N ( {\mathcal{M}_{2}(x_t,u_t,\theta) - \tilde y_t})_2^2 + \rho(\theta) \\
    &\quad  {\sto} \quad  x_{t+1} = \mathcal{M}_{1}(x_t,u_t,\theta), \quad t = 1,\dots,N-1.
\end{aligned} \end{equation}

Problem \eqref{opt_nnss_train} involves $n(N-1)$ optimization variables. Since, for many practical cases $p < n$, i.e., there are fewer outputs than states, this value is greater than the NIO counterpart $p(N-n)$. For this reason, input-output models are computationally preferable when using a constrained optimization approach. Furthermore, as shown in the numerical example of Section~\ref{sec:exampl_mimo}, input-output models are less prone to overfitting than state-space models.

\section{Numerical examples}
\label{sec:exampl}
This section presents numerical examples demonstrating the effectiveness of the proposed method. The examples compare the proposed approach with established alternatives on widely used SI benchmarks. MATLAB code and generated data for all examples are available on the GitHub page: \href{https://github.com/SimonePirrera97/SEM_SysID_with_FLCMO}{\textrm{https://github.com/SimonePirrera97/
SEM\_SysID\_with\_FLCMO}}.

In the following examples, we consider both gray-box and black-box identification problems. For the black-box case, we employ the well-established approach of using neural networks, owing to their ability to approximate any continuous function. Specifically, in all black-box examples, we adopt neural NIO models, i.e., models of the form~\eqref{eq:nio_model} where $\mathcal{M}$ is implemented by a neural network. Following the notation introduced by \cite{sjo94}, we denote the neural NIO model as a nonlinear neural output error (NNOE) model when its training is performed through minimization of the simulation error, and as a nonlinear neural auto-regressive exogenous input (NNARX) model when it is trained by prediction-error minimization. 

In the first example, we compare different neural NIO approaches. The second example validates the effectiveness of the proposed algorithm on a well-established system identification benchmark. The third example compares NNOE training using the proposed method with identification based on recurrent neural networks. Finally, the fourth example addresses gray-box identification of a magnetic levitation system.

The validation metrics adopted here include the root-mean-squared error (RMSE) and the best-fit rate (BFR), defined as:
\begin{equation}
    \text{RMSE} \doteq \sqrt{\frac{1}{ N}\norm{\hat y- y}_2^2}, ~~~ \text{BFR} \doteq \left(1-\frac{\norm{y-\hat{y}}}{\norm{y-m_y}}\right), 
\end{equation}
where $y$ denotes the measured validation output, $\hat y$ the simulated output, and $m_y$ the sample mean of $y$.

Each of the following examples will consider a different value for the parameter $K$, to highlight the algorithm's robustness to its selection. The value of $\tau$ is selected via a trial-and-error procedure to ensure stability of the learning dynamics.

\subsection{Fluid Damper Benchmark}
\label{sec:es01}

This subsection examines the magneto-rheological fluid damper problem from \cite{wang2009identification}. The problem considers magneto-rheological fluid dampers, which are devices used to suppress mechanical vibrations. These dampers contain fluids whose viscosity varies in response to a control signal. The benchmark dataset was obtained by attaching the test damper to a shaker table that generates realistic vibrations. Experimental data from these tests, comprising 2000 data points for training and 1499 data points for testing, are available in MATLAB's System Identification Toolbox.

This subsection evaluates the proposed FL-CMO-based algorithm against leading alternatives for neural NIO modeling in SI, namely NNARX identification via backpropagation and NNOE training with the Levenberg-Marquardt algorithm.

We select a single-layer network with a dynamical order $n=4$ and $6$ neurons. The identification process involved the following methods:
\begin{enumerate}
    \item NNARX: We trained the network using Adam stochastic gradient descent with a batch size of 32. The training algorithm was implemented in Python using TensorFlow’s Keras package, see \cite{chollet2015keras}, with the network trained for 2000 epochs.
    \item NNOE {with Levenberg-Marquardt}: We used the \textit{nnoe} function from the NNSYD MATLAB toolbox, see \cite{norgaad96nnsyd}, which applies the Levenberg-Marquardt algorithm. The gradient and criterion-improvement tolerances were set to $10^{-12}$, while the maximum number of iterations was set to $10^5$.
    \item NNOE with the proposed FL-CMO algorithm: We set the parameters to $\tau = 2 \cdot 10^{-3}$ and $K=1$, with a maximum of $1000$ iterations and regularization $\rho(\theta) = 10^{-3} \norm{\theta}_2^2$.
\end{enumerate}

\begin{table*}[h!]
\caption{{Example 1. Comparison of NNARX (Adam) and NNOE (FL-CMO and Levenberg-Marquardt) identification methods.}}
\label{tab:es_fluid}
	\centering
		\begin{tabular}{|c|cc|}
        \hline
			& \textbf{training time} & \textbf{BFR}\\
   \hline
			{NNOE using FL-CMO} & {106.68 (1.17)} & {88.28 (1.29)} \\
            \hline
			{NNOE using Levenberg-Marquardt)} & 0.32 (0.097) & 66.66 (27.14) \\
			NNARX using Adam & 176.22 ( 32.27) & 68.14 (23.32)\\
            \hline
   \end{tabular}
\end{table*}

We repeat the identification process $10$ times, with different, normally distributed, initializations of the optimization variables. Table \ref{tab:es_fluid} reports the mean value and standard deviation of the training times (measured in seconds) and validation BFRs across all runs for the different algorithms.
The results indicate the FL-CMO algorithm requires lower computational effort compared to Adam while delivering superior performance. In contrast, the Levenberg-Marquardt algorithm significantly reduces computational effort; however, it is less reliable due to the vanishing gradient issue, resulting in poorer average performance.

The obtained results are consistent with or improve upon those reported in Table I of \cite{bemporad22_ekf}, which are based on the same dataset. Specifically, standard BPTT training with the Adam optimizer achieves a validation BFR of 85.5\% for both RNN and LSTM models, while the extended Kalman filter approach proposed in \cite{bemporad22_ekf} yields BFRs of 89.78\% and 88.97\% for RNN and LSTM, respectively.

\subsection{Bouc-Wen system benchmark}
This section considers the Bouc-Wen model, which is described by the dynamics:
\begin{equation}\begin{aligned}
    &m_L \ddot{y}(t) + c_L \dot y(t) + k_L y(t) + z(t) = u(t)\\
    & \dot z(t) = \alpha \dot y(t) -\beta( \gamma \vert  \dot y(t) \vert \vert z(t) \vert^{\nu-1} z(t) + \delta  \dot y(t) \vert z(t) \vert^{\nu} ),
\end{aligned}\end{equation}
where $m_L,k_L,c_L,\alpha,\beta,\gamma,\nu,\delta$ are real scalars. This model is extensively used to characterize mechanical hysteretic effects. The problem of performing its identification belongs to the benchmarks considered during the annual \href{https://www.nonlinearbenchmark.org/benchmarks}{Nonlinear System Identification Benchmarks workshop}. For a comprehensive description of the problem, refer to \cite{noel2016hysteretic}.

A total of 5000 data pairs for training and 5000 data pairs for validation were generated using a multisine input signal, adhering to the provided MATLAB function with its default settings. Training and validation output data were corrupted by additive Gaussian noise with a variance of $\SI{8e-3}{\milli\meter^2}$. For testing, the data provided in the files "{\rm uval\_multisine.mat}" and "{\rm yval\_multisine.mat}" were used, which include 8192 data pairs. 

We train NNOE networks with a dynamical order $n=3$ and two layers, varying the number of neurons from 5 to 10. Training is performed using Algorithm \ref{alg:cap}, with $\tau = 2 \cdot 10^{-3}$, $K=10$, $\epsilon_f,\epsilon_h=10^{-4}$ and regularization $\rho(\theta) = 10^{-3} \norm{\theta}_2^2$. Among the trained networks, the model containing 8 neurons was selected, as it achieved the highest validation accuracy with a $\rm{BFR}$ index of $87.38\%$. The network has 145 parameters, and its training requires $3000$ iterations of Algorithm \ref{alg:cap}, corresponding to $\SI{3521}{\second}$ of computations. The achieved $\rm{BFR}$ on the provided test set is $95.01\%$, demonstrating the model’s generalization capabilities.

\begin{table}[ht]
    \caption{Example 2: Comparison of several methods on the Bouc-Wen system benchmark.}
    \label{tab:esBoucWen}
    \centering
    \begin{tabular}{|c|c|}
    \hline
        \textbf{Method} & \textbf{RMSE} \\
        \hline
        \cite{esfahani2017polynomial} & \SI{1.87e-5}{\milli\meter} \\
        \cite{belz2017automatic} & \SI{16.3e-5}{\milli\meter} \\
        \cite{forgione2021dynonet} & \SI{4.52e-5}{\milli\meter}\\
        \hline
        {NNOE trained by FL-CMO} & {\SI{3.33e-5}{\milli\meter}}\\
        \hline
    \end{tabular}
\end{table}

{Table \ref{tab:esBoucWen} presents the validation RMSE of the trained model and compares it with previously reported results on this benchmark. }
{It is noteworthy that the obtained NNOE model demonstrates superior accuracy compared to the method by \cite{belz2017automatic} and {accuracy comparable to those by \cite{forgione2021dynonet}, using only 5000 training data points rather than 40960, emphasizing the data efficiency of the NNOE model. This performance can be attributed to the NNOE model’s reduced parameterization, which enables effective learning from smaller datasets and mitigates the risk of convergence to poor local minima.} In contrast, accuracy remains lower compared that obtained with the technique introduced by \cite{esfahani2017polynomial}, which is designed specifically for hysteretic systems.}

\subsection{MIMO Wiener-Hammerstein system}
\label{sec:exampl_mimo} 
\begin{figure}[ht]
    \centering
    \includegraphics[width=\linewidth]{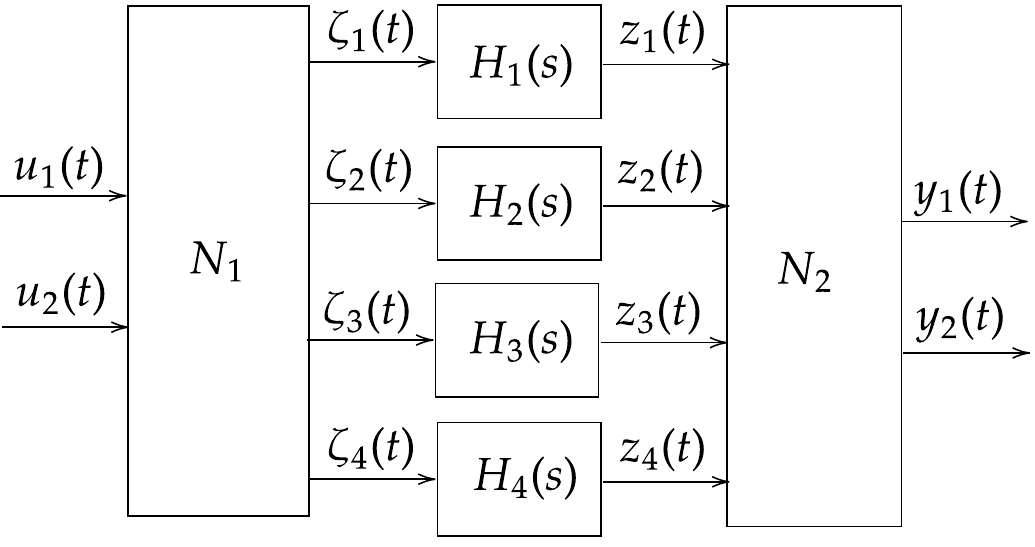}
    \caption{MIMO Wiener-Hammerstein system diagram.}
    \label{fig:whsysdiag}
\end{figure}
\begin{figure*}
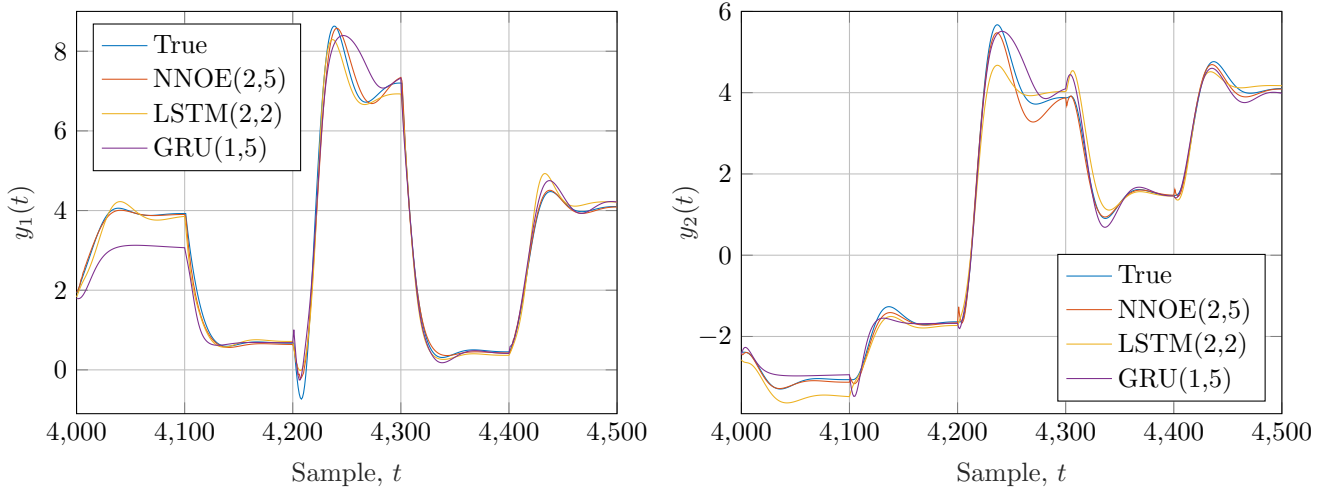

    \input{out1}
    \hfill
    \input{out2}
    \caption{Example 3. Validation output of the top-performing models in each class for output $y_1$ (left) and $y_2$ (right).}
    \label{fig:wh_example_test}
\end{figure*}

In this subsection, we address the identification of the Wiener–Hammerstein system shown in Figure \ref{fig:whsysdiag}, which features two inputs and two outputs. The problem is challenging because the overall behavior is governed by the interplay between linear dynamic blocks and nonlinear components, leading to complex system responses. The input nonlinear block is described by 
\begin{subequations}\label{eq:model_es02a} \begin{align}
\zeta_1 &= 0.2 u_1 u_2 - u_1 - u_2 - 0.04 u_2^2 (u_2+u_1), \\
\zeta_2 &= u_1+u_2, \quad   \zeta_3 = 2 u_2 - 1.3 \zeta_1, \quad \zeta_4 = u_1,
\end{align}
the linear blocks by 
\label{eq:model_es02b} \begin{align}
    H_1(s) &= \frac{16}{s^2 + 40 s + 130},\\
    H_2(s) &= \frac{80}{s^3 + 31 s^2 + 340 s + 2730}, \\
    H_3(s) &= \frac{30}{s + 10}, \quad H_4(s) = \frac{-10}{s^2 + 6 s + 80},
\end{align}
and output nonlinear block by
 \begin{align}\label{eq:model_es02c}
    y_1 &= 7.5 z_4 z_1- 51 z_2+1.02 z_4^2 z_2 - 9 z_1^3 \\ y_2 &=  2.7 z_2 z_3 - 0.729 z_4  z_2^2.
\end{align}\end{subequations}
{We generate a training dataset by simulating the system for $35$ seconds and sampling the obtained signals with period $T_s = \SI{0.01}{\second}$, yielding $N=3500$ training samples. Similarly, we generate a test dataset of $5000$ samples. In both cases, the input signal is a staircase signal with stair length $T_\textrm{stair} = 100 T_s$ and normally distributed amplitudes with unitary variance. 

We train several NNOE models with dynamical order $n=3$ and different neural network structures using the proposed Algorithm \ref{alg:cap} with parameters $\tau = 0.01, K = 100$, $\epsilon_f = \epsilon_h = 10^{-3}$, and the maximum of $500$ iterations.} 
We compare the results obtained using our method with those provided by training LSTM, GRU, and NNARX networks with various architectures. Training of LSTM and GRU is performed using the BPTT algorithm, while the NNARX networks are trained using standard backpropagation (Adam algorithm). To assess model quality, we evaluate the test RMSE and BFR. Table \ref{tab:wh_example_full} summarizes the numerical results concerning training time and accuracy, and marks in green the top-performing models within each class: NNOE(2,5), with two layers of five neurons each; LSTM(2,2), with two layers of two cells; and GRU(1,5), including a single layer of five cells. Figure \ref{fig:wh_example_test} presents the validation outputs of these models. 

{The results reported in Table~\ref{tab:wh_example_full} indicate that the NNARX, LSTM, and GRU models require longer training times and yield lower prediction accuracy than the NNOE model. In particular, the simulation performance of NNARX is significantly worse. This can be attributed to NNARX training relying on PEM, whereas the other models employ SEM. At the same time, the superior performance of NNOE with respect to LSTM and GRU stems from its higher data efficiency and the resilience of its training process to vanishing gradient issues.}

% Table generated by Excel2LaTeX from sheet 'Sheet1'
\begin{table*}
  \centering
  \caption{Example 3: Comparison of several networks.}
    \resizebox{\textwidth}{!}{\begin{tabular}{|cccc|ccccc|cccc|}
    \hline
    %\toprule
    \multicolumn{4}{|c|}{\multirow{2}[2]{*}{\textbf{NNOE network}}} & \multicolumn{5}{c|}{\multirow{2}[2]{*}{\textbf{NNOE Training: FL-CMO algorithm (CPU)}}} & \multicolumn{4}{c|}{\multirow{2}[2]{*}{\textbf{NNOE Validation}}} \\
    \multicolumn{4}{|c|}{}        & \multicolumn{5}{c|}{}                 & \multicolumn{4}{c|}{} \\
    \hline
    \textbf{order} & \textbf{layers} & \textbf{neuron/layer} & \textbf{\# param.} & \textbf{RMSE $y_1$} & \textbf{RMSE $y_2$} & \textbf{iterations} & \textbf{time } & \textbf{time/iter} & \textbf{RMSE $y_1$} & \textbf{RMSE $y_2$} & \textbf{BFR $y_1$} & \textbf{BFR $y_2$} \\
    \hline
    3     & 1     & 5     & 87    & 0.4667 & 0.3141 & 320   & 514.19 & 1.607 & 0.2816 & 0.2670 & 86.03\% & 83.00\% \\
    3     & 1     & 8     & 138   & 0.4411 & 0.2979 & 285   & 633.59 & 2.223 & 0.4557 & 0.3528 & 77.39\% & 77.54\% \\
    3     & 1     & 10    & 172   & 0.3535 & 0.0295 & 500   & 817.81 & 1.636 & 0.1133 & 0.0864 & 94.38\% & 94.50\% \\
    3     & 1     & 15    & 257   & 0.3532 & 0.0252 & 420   & 703.28 & 1.674 & 0.1145 & 0.0925 & 94.32\% & 94.11\% \\
    \hline
    3     & 2     & 3     & 65    & 0.3880 & 0.1733 & 300   & 474.93 & 1.583 & 0.1622 & 0.2078 & 91.95\% & 86.77\% \\
    \rowcolor{green}
    {{3}} & {2} & {5} & {117} & {0.3550} & {0.0728} & {360} & {585.42} & {1.626} & {0.0994} & {0.1071} & {95.07\%} & {93.18\%} \\
    3     & 2     & 7     & 177   & 0.4979 & 0.2764 & 480   & 777.75 & 1.620 & 0.2942 & 0.1718 & 85.41\% & 89.06\% \\
    3     & 2     & 10    & 282   & 0.3541 & 0.0255 & 500   & 826.19 & 1.652 & 0.1188 & 0.1097 & 94.11\% & 93.02\% \\
    \hline
    3     & 3     & 3     & 77    & 0.4049 & 0.1898 & 260   & 453.32 & 1.744 & 0.1494 & 0.1958 & 92.59\% & 87.54\% \\
    3     & 3     & 5     & 147   & 0.4247 & 0.2683 & 180   & 438.12 & 2.434 & 0.2694 & 0.2975 & 86.64\% & 81.06\% \\
    3     & 3     & 7     & 233   & 0.6969 & 0.7338 & 140   & 244.47 & 1.746 & 0.2452 & 0.2378 & 87.84\% & 84.86\% \\
    \hline
    \multicolumn{1}{r}{} &       &       & \multicolumn{1}{r}{} &       &       &       &       & \multicolumn{1}{r}{} &       &       &       &  \\
    \hline
    \multicolumn{4}{|c|}{\multirow{2}[2]{*}{\textbf{LSTM network} }} & \multicolumn{5}{c|}{\multirow{2}[2]{*}{\textbf{LSTM Training: Adam BPTT algorithm (GPU)}}} & \multicolumn{4}{c|}{\multirow{2}[2]{*}{\textbf{LSTM Validation}}} \\
    \multicolumn{4}{|c|}{}        & \multicolumn{5}{c|}{}                 & \multicolumn{4}{c|}{} \\
    \hline
    \textbf{order} & \textbf{layers} & \textbf{units/layer} & \textbf{\# param.} & \textbf{RMSE $y_1$} & \textbf{RMSE $y_2$} & \textbf{epochs} & \textbf{time } & \textbf{time/epoch} & \textbf{RMSE $y_1$} & \textbf{RMSE $y_2$} & \textbf{BFR $y_1$} & \textbf{BFR $y_2$} \\
    \hline
    4     & 1     & 2     & 46    & 0.5475 & 0.4164 & 12000 & 841.05 & 0.0701 & 0.6236 & 0.7680 & 69.06\% & 51.11\% \\
    6     & 1     & 3     & 80    & 0.4142 & 0.3207 & 12000 & 821.17 & 0.0684 & 0.5876 & 0.6221 & 70.85\% & 60.40\% \\
    8     & 1     & 4     & 122   & 0.2371 & 0.1841 & 12000 & 1036.00 & 0.0863 & 0.4236 & 0.2769 & 78.99\% & 82.37\% \\
    10    & 1     & 5     & 172   & 0.0793 & 0.0790 & 12000 & 827.54 & 0.0690 & 0.3811 & 0.2150 & 81.09\% & 86.31\% \\
    14    & 1     & 7     & 296   & 0.0545 & 0.0505 & 12000 & 826.17 & 0.0688 & 0.3606 & 0.2628 & 82.11\% & 83.27\% \\
    \hline
    \rowcolor{green} 
    {8} & {2} & {2} & {86} & {0.2747} & {0.2142} & {12000} & {1690.19} & {0.1408} & {0.3254} & {0.2512} & {83.86\%} & {84.00\%} \\
    12    & 2     & 3     & 164   & 0.2131 & 0.1636 & 12000 & 3061.88 & 0.2552 & 0.3858 & 0.4607 & 80.86\% & 70.67\% \\
    16    & 2     & 4     & 266   & 0.1781 & 0.1543 & 12000 & 2879.37 & 0.2399 & 0.8467 & 0.4391 & 58.00\% & 72.04\% \\
    20    & 2     & 5     & 392   & 0.1867 & 0.1390 & 12000 & 2392.98 & 0.1994 & 0.8459 & 0.7204 & 58.04\% & 54.14\% \\
    \hline
    12    & 3     & 2     & 126   & 0.2046 & 0.2049 & 12000 & 2865.79 & 0.2388 & 1.7057 & 0.3697 & 15.38\% & 76.46\% \\
    18    & 3     & 3     & 248   & 0.0755 & 0.0818 & 12000 & 2309.90 & 0.1925 & 0.3655 & 0.1999 & 81.87\% & 87.27\% \\
    \hline
    \multicolumn{1}{c}{} &       &       & \multicolumn{1}{c}{} &       &       &       &       & \multicolumn{1}{c}{} &       &       &       & \multicolumn{1}{c}{} \\
    \hline
    \multicolumn{4}{|c|}{\multirow{2}[2]{*}{\textbf{GRU network} }} & \multicolumn{5}{c|}{\multirow{2}[2]{*}{\textbf{GRU Training: Adam BPTT algorithm (GPU)}}} & \multicolumn{4}{c|}{\multirow{2}[2]{*}{\textbf{GRU Validation}}} \\
    \multicolumn{4}{|c|}{}        & \multicolumn{5}{c|}{}                 & \multicolumn{4}{c|}{} \\
    \hline
    \textbf{order} & \textbf{layers} & \textbf{units/layer} & \textbf{\# param.} & \textbf{RMSE $y_1$} & \textbf{RMSE $y_2$} & \textbf{epochs} & \textbf{time } & \textbf{time/epoch} & \textbf{RMSE $y_1$} & \textbf{RMSE $y_2$} & \textbf{BFR $y_1$} & \textbf{BFR $y_2$} \\
    \hline
    3     & 1     & 3     & 71    & 0.4960 & 0.3020 & 12000 & 1498.09 & 0.1248 & 1.1727 & 1.3030 & 41.82\% & 17.05\% \\
    4     & 1     & 4     & 106   & 0.1376 & 0.1128 & 12000 & 1297.46 & 0.1081 & 0.3307 & 0.2031 & 83.59\% & 87.07\% \\
    \rowcolor{green} 
    {5} & {1} & {5} & {147} & {0.1617} & {0.1380} & {12000} & {819.22} & {0.0683} & {0.1929} & {0.2127} & {90.43\%} & {86.46\%} \\
    7     & 1     & 7     & 247   & 0.0564 & 0.0504 & 12000 & 810.15 & 0.0675 & 0.2219 & 0.2425 & 88.99\% & 84.56\% \\
    10    & 1     & 10    & 442   & 0.0741 & 0.0853 & 12000 & 1460.55 & 0.1217 & 0.4643 & 0.7117 & 76.96\% & 54.69\% \\
    \hline
    6     & 2     & 3     & 143   & 0.2337 & 0.1584 & 12000 & 1578.26 & 0.1315 & 0.3226 & 0.3752 & 84.00\% & 76.11\% \\
    8     & 2     & 4     & 226   & 0.1353 & 0.1265 & 12000 & 1581.26 & 0.1318 & 0.3751 & 0.5469 & 81.39\% & 65.18\% \\
    10    & 2     & 5     & 327   & 0.0326 & 0.0309 & 12000 & 1591.46 & 0.1326 & 0.6165 & 0.2345 & 69.41\% & 85.07\% \\
    14    & 2     & 7     & 583   & 0.0253 & 0.0325 & 12000 & 1571.18 & 0.1309 & 0.3244 & 0.3920 & 83.91\% & 75.05\% \\
    \hline
    6     & 3     & 2     & 114   & 0.1929 & 0.1647 & 12000 & 3638.46 & 0.3032 & 0.3943 & 0.2501 & 80.44\% & 84.08\% \\
    9     & 3     & 3     & 215   & 0.0829 & 0.0858 & 12000 & 3741.53 & 0.3118 & 0.3664 & 0.2032 & 81.83\% & 87.07\% \\
    \hline
    \multicolumn{1}{c}{} &       &       & \multicolumn{1}{c}{} &       &       &       &       & \multicolumn{1}{c}{} &       &       &       & \multicolumn{1}{c}{} \\
    \hline
    \multicolumn{4}{|c|}{\multirow{2}[2]{*}{\textbf{NNARX network}}} & \multicolumn{5}{c|}{\multirow{2}[2]{*}{\textbf{NNARX Training: Adam BP algorithm (GPU)} }} & \multicolumn{4}{c|}{\multirow{2}[2]{*}{\textbf{NNARX Validation}}} \\
    \multicolumn{4}{|c|}{}        & \multicolumn{5}{c|}{}                 & \multicolumn{4}{c|}{} \\
    \hline
    \textbf{order} & \textbf{layers} & \textbf{neuron/layer} & \textbf{\# param.} & \textbf{RMSE $y_1$} & \textbf{RMSE $y_2$} & \textbf{epochs} & \textbf{time } & \textbf{time/epoch} & \textbf{RMSE $y_1$} & \textbf{RMSE $y_2$} & \textbf{BFR $y_1$} & \textbf{BFR $y_2$} \\
    \hline
    3     & 1     & 5     & 87    & 5.0156 & 4.0211 & 1000  & 89.59 & 0.090 & 6.3442 & 3.6032 & -214.72\% & -129.39\% \\
    3     & 1     & 8     & 138   & 6.8757 & 6.9953 & 1000  & 93.47 & 0.093 & 7.2585 & 6.9879 & -260.08\% & -344.87\% \\
    3     & 1     & 10    & 172   & 7.8516 & 7.7647 & 1000  & 104.78 & 0.105 & 7.7406 & 9.3001 & -284.00\% & -492.07\% \\
    3     & 1     & 15    & 257   & 8.0630 & 9.8919 & 1000  & 112.59 & 0.113 & 7.4535 & 10.7121 & -269.76\% & -581.97\% \\
    \hline
    3     & 2     & 3     & 65    & 8.2622 & 9.5509 & 1000  & 170.33 & 0.170 & 8.1210 & 9.6828 & -302.87\% & -516.44\% \\
    3     & 2     & 5     & 117   & 7.7178 & 9.4936 & 1000  & 120.02 & 0.120 & 6.6021 & 9.3484 & -227.52\% & -495.15\% \\
    3     & 2     & 7     & 177   & 9.9773 & 8.4243 & 1000  & 117.51 & 0.118 & 9.8954 & 8.5393 & -390.90\% & -443.64\% \\
    3     & 2     & 10    & 282   & 2.5683 & 5.5859 & 1000  & 105.30 & 0.105 & 2.0974 & 5.8676 & -4.05\% & -273.55\% \\
    \hline
    3     & 3     & 3     & 77    & 5.4776 & 3.1718 & 1000  & 120.97 & 0.121 & 5.5254 & 2.8112 & -174.11\% & -78.97\% \\
    3     & 3     & 5     & 147   & 5.3081 & 5.7561 & 1000  & 112.76 & 0.113 & 5.2680 & 6.5923 & -161.34\% & -319.69\% \\
    3     & 3     & 7     & 233   & 7.6568 & 9.9372 & 1000  & 120.09 & 0.120 & 7.5069 & 11.4208 & -272.40\% & -627.09\% \\
    \bottomrule
    \end{tabular}}%
  \label{tab:wh_example_full}%
\end{table*}%

Concerning computational effort, the columns time/iter and time/epoch in Table~\ref{tab:wh_example_full} report the average computation time per iteration of Algorithm~\ref{alg:cap} and per epoch of BPTT, respectively. From these results, we observe that a single iteration of Algorithm~\ref{alg:cap} is more computationally demanding than one BPTT epoch for LSTM and GRU training. However, training NNOEs with Algorithm~\ref{alg:cap} requires significantly fewer iterations overall, resulting in a lower total training time. This improvement is primarily due to the absence of vanishing gradient issues.

\subsection{Magnetic levitation system gray-box identification}
We analyze the magnetic levitation system described by the continuous-time model
\begin{equation}\begin{aligned}
 \ddot z &= g - \frac{k_m i^2 + k_0}{m z^2},
\end{aligned} \label{eq:levit} \end{equation}
where $z$ denotes the distance between the magnet and the ball, $m$ is the mass of the ball, $g$ denotes the gravitational acceleration, and $k_m$ and $k_0$ are constants determined by the magnet's characteristics. The objective of the identification procedure is to estimate $k_m$ and $k_0$, with $m=\SI{24.197}{\gram}$ and $g=\SI{9.81}{\meter\per\second^2}$ treated as known constants. 

We simulate the differential equation \eqref{eq:levit} and collect $N = 200$ input-output data points using a sampling rate of $T_s = \SI{10}{\milli\second}$. The output data is subsequently corrupted by additive measurement noise uniformly distributed within the interval %$[\SI{-5}{\milli\meter},\SI{+5}{\milli\meter}]$.
$[-5,5] $ mm.

We discretize the dynamic system equation using Euler’s method, resulting in the following model
\begin{equation}\label{eq:dt_mdl_levit}
 m\, z_{t-2}^2 \left( \frac{ z_{t} - 2 z_{t-1} + z_{t-2}}{T_s^2} - g\right) + k_m\, i_{t-2}^2 + k_0 = 0.
\end{equation}
According to the proposed approach, we evaluate the Jacobian of the constraints defined in Equation \eqref{eq:dt_mdl_levit}. Specifically, the computation requires determining the partial derivatives:
\begin{subequations}
\begin{align}
    \frac{\partial h_t}{\partial k_m} &= i_{t-2}^2, \qquad 
    \frac{\partial h_t}{\partial z_{t}} = \frac{m}{T_s^2} z_{t-2}^2, \\
    \frac{\partial h_t}{\partial z_{t-1}} &= -\frac{2m}{T_s^2} z_{t-2}^2,\qquad\frac{\partial h_t}{\partial k_0} = 1, \\
    \frac{\partial h_t}{\partial z_{t-2}} &= \frac{m}{T_s^2} z_{t-2} (2z_t -4z_{t-1} + 3z_{t-2} ) - 2mg z_{t-2}.
\end{align}
\end{subequations}
{We apply Algorithm \ref{alg:cap} using parameters $K = 2$ and $\tau = 10^{-3}$. }For comparison, the system is also identified using the following methods:
\begin{enumerate}
    \item Least-Squares (LS): The one-step-ahead prediction-error was minimized, resulting in a linear and convex LS problem. The LS approach is widely regarded as the standard method for estimating the physical parameters of dynamical systems, particularly when the system equation depends linearly on the unknown parameters.
    \item Adam {Backpropagation Through Time (BPTT)}: The unconstrained SEM optimization problem, as outlined in Equation \eqref{eq:uncons_formulation}, was formulated and solved using the Adam algorithm. The loss gradient was computed through BPTT, using a learning rate of $10^{-3}$, and optimization was conducted over $10^4$ epochs.
\end{enumerate}

\begin{table}[ht]
    \centering
    \caption{{Example 4: Comparison of parameters estimated for the magnetic levitation system using SEM with FL-CMO, SEM with Adam BPTT, and PEM based on LS.}}
    \label{tab:ese_levit_results}
    \begin{tabular}{|p{2.5cm}|c c|}
    \hline
    & $k_m$ (\rm{H\,m}) & $k_0$ (\rm{N\,m}) \\
    \hline
    {True} & $2.1039\times 10^{-4}$ & $0.0$ \\
    \hline
    {FL-CMO} & ${2.1037\times 10^{-4}}$ & ${2.4108\times 10^{-5}}$ \\
    {Adam BPTT}  & ${7.46\times 10^{-4}}$  & ${-0.143}$\\
    {LS} & ${2.1228\times 10^{-4}}$ & ${-3.6130\times 10^{-4}}$ \\
    \hline
    \end{tabular}
\end{table}

We repeat the optimization $10$ times for each method, using distinct, zero-mean, normally distributed initializations for the optimization variables. In Table \ref{tab:ese_levit_results}, we compare the results of all three methods, selecting the best estimate among all runs for each approach. The results indicate that FL-CMO demonstrates the highest identification accuracy. Although the LS method yields an acceptable estimate, the accuracy of the $k_0$ parameter is one order of magnitude lower than that of the FL-CMO approach. Conversely, the Adam BPTT method produces inaccurate estimates for both parameters. We explain these results by observing that, although both FL-CMO and Adam BPTT are based on SEM, Adam BPTT suffers from the vanishing gradient issue.

\section{Conclusions}
\label{sec:concl}
This study investigates system identification for nonlinear input-output models by formulating the simulation-error minimization problem as a constrained optimization task. This formulation avoids vanishing and exploding gradients, which are common issues that hinder convergence in traditional approaches. To solve the constrained problem, we introduce a novel identification algorithm that integrates feedback linearization controlled multipliers optimization with the Euler discretization method, enabling efficient and accurate model estimation.

Our theoretical analysis provides guarantees of local convergence and establishes connections between the proposed approach and the conventional unconstrained formulation. By leveraging the problem structure and employing sparse Q-less QR factorization, we reduce the computational complexity of the iterative process compared to standard unconstrained optimization techniques. 

Extensive numerical experiments on black-box and gray-box problems reveal that the proposed approach outperforms conventional methods. Additionally, we find that adequately trained neural input-output models can exhibit superior performance compared to state-space neural ones, such as LSTM and GRU. These findings offer an efficient alternative to gradient-based identification methods for nonlinear systems. Future research will focus on developing solutions to further enhance efficiency in large-scale problems.

\bibliographystyle{apalike}
\bibliography{biblio}
\end{document}